\documentclass[11pt,a4paper]{article}

\usepackage[utf8]{inputenc}
\usepackage[T1]{fontenc}
\usepackage{lmodern}
\usepackage{amsmath,amssymb,amsthm}
\usepackage{graphicx}
\usepackage[margin=1in]{geometry}
\usepackage{microtype}
\usepackage[hidelinks]{hyperref}
\usepackage{authblk}
\usepackage{titling}

\title{A Scaling-Parameter Framework for Perimeter and Area\\
in Self-Similar Planar Fractals}
\author{Pedro Marotta\thanks{Correspondence: \href{mailto:pedromarottaar@gmail.com}{pedromarottaar@gmail.com}.}}
\affil{St.~Andrews Scots School, Olivos, Vicente Lopez Partido,\protect\\ Buenos Aires Province CP~1636, Argentina}
\date{}

\begin{document}

\maketitle

\noindent\textit{Author's accepted manuscript. Published in: Marotta, P.\ (2026).
A Scaling-Parameter Framework for Perimeter and Area in Self-Similar Planar
Fractals. \emph{Journal of High School Science}, 10(2), 275--295.
DOI: \href{https://doi.org/10.64336/001c.162173}{10.64336/001c.162173}.
Licensed under CC~BY-NC-SA~4.0.}

\bigskip

\begin{abstract}
\noindent The Koch snowflake is a classical example of a planar curve with
infinite perimeter enclosing a finite, positive area. Although such examples
are well known individually, classical treatments typically analyze each
construction in isolation and classify them by similarity dimension. This
paper develops a unified parameter-space representation for a class of
self-similar planar constructions, organized by two integers~--- the number
of self-similar pieces $N$ and the inverse linear scale factor $r$~--- together
with two derived growth ratios $\alpha = N/r$ and $\beta = N/r^{2}$, which
govern perimeter scaling and area scaling respectively. The condition
$1 < D = \log N / \log r < 2$ for non-integer planar fractal dimension is
classical; the present framework recasts these results in coordinates that
make perimeter behavior and area behavior directly comparable observables
within a single coordinate system. The $(N, r)$ parameter space is
partitioned into three regimes~--- $N \le r$, $r < N < r^{2}$, and
$N \ge r^{2}$~--- corresponding to qualitatively distinct asymptotic behaviors
of perimeter and area jointly. The framework is then refined by a
construction-class distinction: within the same intermediate regime
$r < N < r^{2}$, additive constructions (in which the object of interest is
the region bounded by the iterated curve) yield positive finite asymptotic
area under a stated non-overlap assumption, while subtractive constructions
(in which the object is the iterated set itself) yield zero asymptotic area.
This refinement records a structural non-equivalence inside the same
dimension class that is not visible from $D$ alone. Four worked examples
illustrate the framework~--- the Sierpinski triangle, the Sierpinski carpet,
the Koch snowflake, and a Koch-style construction on a square invented by the
author~--- and four further constructions are analyzed predictively to
demonstrate that the framework's diagnostic outputs follow from
$(N, r, \text{construction class})$ without re-derivation: one example in
each of the three regimes, plus a pair with identical $(N, r)$, and therefore
identical similarity dimension, that nevertheless exhibit qualitatively
different asymptotic area outcomes because of differing construction class.
The contribution of the paper lies in formulation and synthesis rather than
in new mathematics: it consolidates several classical results into a single
diagnostic representation in which, given $(N, r)$ and a specification of
construction class, the asymptotic behavior of perimeter and area can be
inferred directly without re-deriving each example.
\end{abstract}


\section{Introduction}

A curve of infinite length that bounds a finite area is one of the more counterintuitive constructions in classical geometry. The Koch snowflake, introduced by Helge von Koch in 1904 \cite{vonKoch1904}, is the canonical example: at each iteration, every line segment of the boundary is replaced by four segments of one-third the length, and the resulting curve in the limit has infinite arc length while the region it encloses has area $2\sqrt{3}/5$ of an initial unit-side equilateral triangle. Comparable behavior is exhibited by other self-similar planar constructions, including the Sierpinski triangle \cite{Sierpinski1915} and the Sierpinski carpet, although the precise sense in which ``perimeter'' is taken to diverge and whether the limiting object has positive or zero area depend on construction details that are not always made explicit in introductory treatments.

The textbook condition under which a planar self-similar fractal exhibits non-integer dimension is well known: the similarity dimension $D = \log N / \log r$ lies between $1$ and $2$ \cite{Falconer2014,Edgar2008}. The four examples listed above all satisfy $1 < D < 2$ and all have iteration-$n$ perimeter that diverges as a function of $n$ under suitable conventions. They do not, however, all share the ``infinite perimeter and finite (positive) area'' property of the Koch snowflake: the Sierpinski triangle and Sierpinski carpet have two-dimensional Lebesgue measure equal to zero when viewed as the iterated set itself.

The aim of this paper is to provide a unified parameter-space representation for these constructions that (i) parameterizes them by two integers $N$ and $r$ and by two derived growth ratios $\alpha = N/r$ and $\beta = N/r^{2}$, (ii) derives perimeter and area scaling explicitly under perimeter conventions stated per construction class, (iii) classifies the resulting behavior into three regimes in $(N, r)$-space and identifies the intermediate-dimension regime $r < N < r^{2}$ (equivalently $1 < D < 2$), and (iv) within that regime, distinguishes between two construction classes~--- additive (where the object of interest is a region bounded by the iterated curve) and subtractive (where the object of interest is the iterated set itself)~--- and records that the asymptotic two-dimensional Lebesgue measure depends on the class even though similarity dimension does not.

The contribution of the paper is one of formulation and synthesis rather than of new mathematics. The dimension formula $D = \log N / \log r$, the condition $1 < D < 2$ for planar fractal dimension, and Hutchinson's theorem on the two-dimensional Lebesgue measure of self-similar sets with $D < 2$ are all classical \cite{Falconer2014,Edgar2008,Hutchinson1981}. What the present framework adds is (a) a change of representation in which $\alpha = N/r$ and $\beta = N/r^{2}$ become the primary control parameters and the regime classification follows from their joint behavior; (b) a synthesis in which four canonical examples are located within a single coordinate system rather than treated in isolation; and (c) a construction-class refinement that records, within the same intermediate regime, a structural non-equivalence between additive and subtractive constructions that is not apparent from similarity dimension alone. The resulting framework is intended to be used as a diagnostic: given $(N, r)$ and a specification of construction class, the asymptotic behavior of perimeter and of area can be read off from $\alpha$, $\beta$, and the construction class without re-deriving each example.

Section~\ref{sec:preliminaries} introduces preliminaries and conventions. Section~\ref{sec:examples} presents the four worked examples. Section~\ref{sec:framework} develops the $(N, r)$ framework, defines the derived ratios $\alpha$ and $\beta$, presents the regime classification, states the construction-class refinement and the diagnostic use of the framework, and applies the framework predictively (\S\ref{sec:predictive}) to four further constructions outside the canonical four of Section~\ref{sec:examples}~--- one in each of the three regimes plus a pair with identical $(N, r)$ and differing construction class. Section~\ref{sec:regimefigure} contains the regime figure. Section~\ref{sec:discussion} discusses limitations and the relation to construction classes outside the present scope, including stochastic fractal models. Section~\ref{sec:conclusion} concludes.

\section{Background}\label{sec:preliminaries}

\subsection{Iterated self-similar planar construction}

A self-similar planar construction in this paper is generated by an iterative replacement rule. At iteration zero, an initial figure is given (a line segment, a polygon, or a region). At each subsequent iteration, every segment of the current figure is replaced by $N$ copies of itself, each scaled by a factor of $1/r$ (with $r > 1$ a positive real), arranged in a fixed geometric pattern relative to the original segment. The four examples treated in Section~\ref{sec:examples} are special cases of this scheme.

\subsection{Construction classes: additive and subtractive}

Two construction classes are distinguished. In an \emph{additive} construction, each line segment of the iterated boundary is replaced by $N$ shorter segments arranged so as to form an outward bump; the object of interest at iteration $n$ is the region of the plane bounded by the iterated curve. The Koch snowflake and the Koch-style construction on a square presented in Section~\ref{sec:koch-square} are additive. In a \emph{subtractive} construction, the figure at iteration zero is a filled two-dimensional region; at each iteration, that region is partitioned into smaller similar pieces and a fixed subset of the pieces is removed. The object of interest is the resulting set of remaining pieces. The Sierpinski triangle and the Sierpinski carpet are subtractive.

\subsection{Similarity dimension}

For a self-similar construction generated by $N$ pieces each scaled by $1/r$, the similarity dimension $D$ is defined by
\begin{equation*}
D = \frac{\log N}{\log r}. \tag{1}
\end{equation*}
For self-similar iterated function systems satisfying the open set condition, $D$ coincides with the Hausdorff dimension of the limit set \cite{Hutchinson1981}. The Hausdorff dimension itself is defined measure-theoretically \cite[Chapter~2]{Falconer2014} and is in general distinct from the similarity dimension; the present paper computes the similarity dimension only and notes the equivalence under the open set condition for the four examples treated.

\subsection{Perimeter convention}\label{sec:perimeter-convention}

The ``perimeter'' of a self-similar construction is convention-dependent and the convention is stated explicitly per construction class in this paper.

For an additive construction, the perimeter at iteration $n$ is taken to be the arc length of the iterated curve. Each iteration multiplies the number of segments by $N$ and divides each segment length by $r$, so
\begin{equation*}
P_n = P_0 \cdot \left(\frac{N}{r}\right)^n, \tag{2}
\end{equation*}
where $P_0$ is the arc length of the initial figure.

For a subtractive construction, the perimeter at iteration $n$ is taken to be the total edge length of all sub-pieces present at iteration $n$. This is a non-standard convention in the broader literature: alternatives include the boundary of the convex hull of the iterated set (which is constant for the examples here) and the Hausdorff one-measure of the limit set (which is infinite for sets of dimension greater than one, so does not yield a finite-iteration sequence). Under the convention adopted here,
\begin{equation*}
P_n = P_0 \cdot \left(\frac{N}{r}\right)^n, \tag{3}
\end{equation*}
where $P_0$ is the total edge length at iteration zero.

The convention adopted for subtractive constructions matches what is computed informally in many introductory presentations as the iterated ``perimeter'' and yields the same scaling form as~(2). The asymmetry~--- that the limit object's Hausdorff one-measure is not the limit of $P_n$ in the subtractive case~--- is acknowledged as a limitation and discussed in Section~\ref{sec:discussion}.

\subsection{Area}

The area at iteration $n$ of an additive construction is
\begin{equation*}
A_n^{\mathrm{add}} = A_0 + C \cdot \sum_{k=1}^{n} \left(\frac{N}{r^{2}}\right)^{k-1}. \tag{4}
\end{equation*}

Under a stated non-overlap assumption (each added bump intersects neither the existing iterated region nor any other added bump beyond its base), the total area added at successive iterations forms a geometric sequence. From one iteration to the next, the number of segments is multiplied by $N$ while each added bump has $1/r^{2}$ the area of those added at the previous iteration. The total added area is therefore multiplied by $N/r^{2}$ at each step, and the area sum reduces to a geometric series with ratio $N/r^{2}$.

The area at iteration $n$ of a subtractive construction is
\begin{equation*}
A_n = A_0 \cdot \left(\frac{N}{r^{2}}\right)^n, \tag{5}
\end{equation*}
since at each iteration the construction retains a fraction $N/r^{2}$ of the previous area ($N$ pieces of area $1/r^{2}$ of the previous each).

In both cases the asymptotic behavior as $n \to \infty$ is governed by the ratio $N/r^{2}$.

\section{Worked examples}\label{sec:examples}

\subsection{Sierpinski triangle}\label{sec:sierpinski-triangle}

\begin{figure}[htbp]
  \centering
  \includegraphics[width=0.95\textwidth]{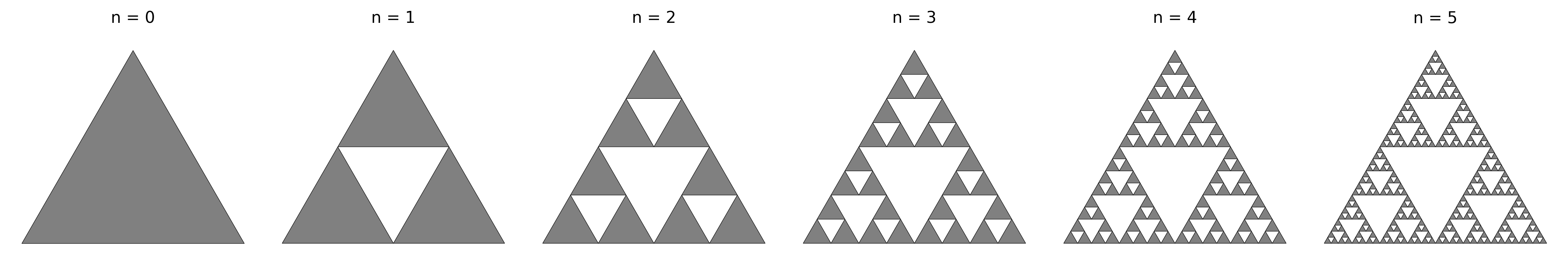}
  \caption{\emph{Sierpinski triangle, iterations $n=0$ through $n=5$.} The construction starts from a closed equilateral triangle ($n=0$). At each iteration, every triangle of side length $s$ is replaced by three triangles of side length $s/2$ placed at its vertices, with the central inverted triangle removed. The number of self-similar pieces is $N=3$, the linear scale factor is $r=2$, and the similarity dimension is $D = \log 3 / \log 2 \approx 1.585$. The set is subtractive: at each step, area is removed and the asymptotic two-dimensional Lebesgue measure of the limit set is zero.}
  \label{fig:sierpinski-triangle}
\end{figure}

The Sierpinski triangle was described geometrically by Wac{\l}aw Sierpi{\'n}ski in 1915 \cite{Sierpinski1915}, although the pattern itself appeared as ornamentation many centuries earlier. Starting from an equilateral triangle of side length~$1$, each iteration partitions every triangle into four congruent equilateral sub-triangles by joining the midpoints of its sides, and removes the central sub-triangle (Figure~\ref{fig:sierpinski-triangle}).

The construction is subtractive with $N = 3$ (three sub-triangles retained per iteration) and $r = 2$ (each retained sub-triangle has side length one-half of its parent). The similarity dimension is
\begin{equation*}
D = \frac{\log 3}{\log 2} \approx 1.585. \tag{6}
\end{equation*}
Under the perimeter convention adopted in \S\ref{sec:perimeter-convention} for subtractive constructions, the total edge length at iteration $n$ is the number of sub-triangles times the perimeter of each, namely
\begin{equation*}
P_n = 3^{n} \cdot (3 \cdot 2^{-n}) = 3 \cdot \left(\frac{3}{2}\right)^{n}, \tag{7}
\end{equation*}
which diverges as $n \to \infty$ since $3/2 > 1$. The boundary of the convex hull of the iterated set, by contrast, is the original triangle's boundary and remains equal to $3$ at every iteration; the divergence claimed here therefore concerns the total edge-length convention adopted, not the convex-hull boundary.

The area at iteration $n$ is
\begin{equation*}
A_n = \frac{\sqrt{3}}{4} \cdot \left(\frac{3}{4}\right)^{n}, \tag{8}
\end{equation*}
since at each iteration three retained sub-triangles each occupy area $1/4$ of the previous triangle's area, so retained area shrinks by a factor $3/4$. As $n \to \infty$, $A_n \to 0$. The Sierpinski triangle as a set therefore has two-dimensional Lebesgue measure zero, consistent with the general result that any self-similar set with similarity dimension $D < 2$ has zero two-dimensional Lebesgue measure \cite{Hutchinson1981}.

\subsection{Sierpinski carpet}\label{sec:sierpinski-carpet}

\begin{figure}[htbp]
  \centering
  \includegraphics[width=0.95\textwidth]{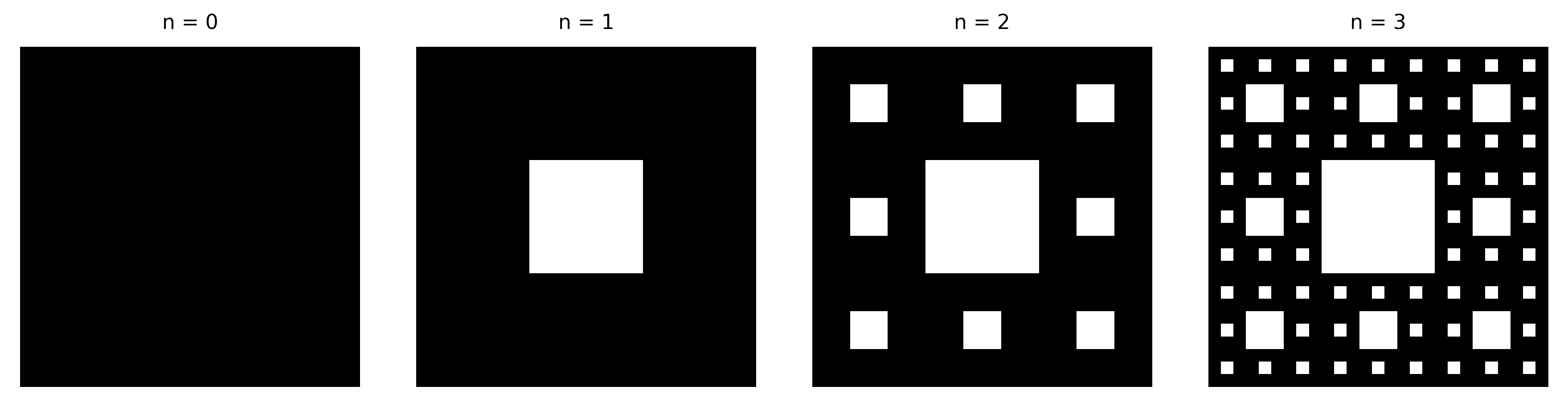}
  \caption{\emph{Sierpinski carpet, iterations $n=0$ through $n=3$.} The construction starts from a closed unit square ($n=0$). At each iteration, every square of side $s$ is partitioned into a $3\times 3$ grid of nine sub-squares of side $s/3$, and the central sub-square is removed. The number of retained pieces is $N=8$, the linear scale factor is $r=3$, and the similarity dimension is $D = \log 8 / \log 3 \approx 1.893$. The set is subtractive; the asymptotic two-dimensional Lebesgue measure of the limit set is zero.}
  \label{fig:sierpinski-carpet}
\end{figure}

The Sierpinski carpet, introduced by Sierpi{\'n}ski in 1916, applies an analogous subtractive rule to a square. The unit square is partitioned into nine congruent sub-squares of side $1/3$; the central sub-square is removed; the rule is applied to each of the eight retained sub-squares, and so on (Figure~\ref{fig:sierpinski-carpet}).

The construction is subtractive with $N = 8$ and $r = 3$. The similarity dimension is
\begin{equation*}
D = \frac{\log 8}{\log 3} \approx 1.893. \tag{9}
\end{equation*}
The total edge length at iteration $n$, under the convention of \S\ref{sec:perimeter-convention}, is
\begin{equation*}
P_n = 4 \cdot \left(\frac{8}{3}\right)^{n}, \tag{10}
\end{equation*}
which diverges. The area at iteration $n$ is
\begin{equation*}
A_n = \left(\frac{8}{9}\right)^{n}, \tag{11}
\end{equation*}
which tends to zero.

\subsection{Koch snowflake}\label{sec:koch-snowflake}

\begin{figure}[htbp]
  \centering
  \includegraphics[width=0.95\textwidth]{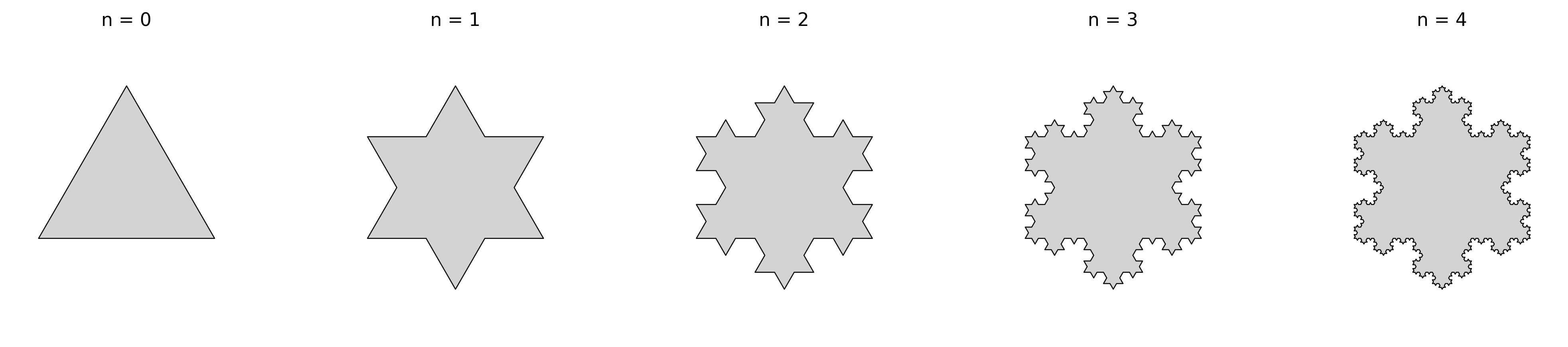}
  \caption{\emph{Koch snowflake, iterations $n=0$ through $n=4$.} The construction starts from a closed equilateral triangle of side length~$1$ ($n=0$). At each iteration, every boundary segment of length $s$ is divided into three equal parts and the middle part is replaced by the two other sides of an outward-pointing equilateral triangle, so that each segment is replaced by four segments of length $s/3$. The generator parameters are $N=4$, $r=3$, and the similarity dimension of the limit curve is $D = \log 4 / \log 3 \approx 1.262$. The construction is additive: the perimeter of the iterated curve diverges as $(4/3)^{n}$, while the area enclosed by the curve converges to $2\sqrt{3}/5 \approx 0.693$ under the standard non-overlap of added bumps.}
  \label{fig:koch-snowflake}
\end{figure}

The Koch snowflake \cite{vonKoch1904} is generated by an additive rule applied to the boundary of an equilateral triangle of side length~$1$. At each iteration, every segment of the current boundary is divided into three equal parts; the middle third is replaced by the other two sides of an equilateral triangle that is built outward on that middle third; the original middle third is then removed (Figure~\ref{fig:koch-snowflake}).

The construction is additive with $N = 4$ (each segment is replaced by four shorter segments) and $r = 3$ (each new segment is one-third the length of its parent). The similarity dimension is
\begin{equation*}
D = \frac{\log 4}{\log 3} \approx 1.262. \tag{12}
\end{equation*}
The arc length at iteration $n$ is
\begin{equation*}
P_n = 3 \cdot \left(\frac{4}{3}\right)^{n}, \tag{13}
\end{equation*}
which diverges since $4/3 > 1$.

The area at iteration $n$ is computed as the initial triangle's area plus the sum of bump areas added at each iteration. The initial area is $\sqrt{3}/4$. At iteration~1, three new equilateral triangles of side $1/3$ are added, each of area $(1/3)^{2} \cdot \sqrt{3}/4$. At iteration~2, every segment of the iteration-1 boundary (twelve segments of length $1/3$) yields one additional outward triangle of side $1/9$, so twelve new triangles of area $(1/9)^{2} \cdot \sqrt{3}/4$ are added. In general, at iteration $k \ge 1$, $3 \cdot 4^{k-1}$ new triangles of side $1/3^{k}$ are added, contributing total area
\begin{equation*}
\text{(area added at iteration } k) = 3 \cdot 4^{k-1} \cdot \left(\frac{1}{3^{k}}\right)^{2} \cdot \frac{\sqrt{3}}{4} = \frac{\sqrt{3}}{12} \cdot \left(\frac{4}{9}\right)^{k-1}. \tag{14}
\end{equation*}
Summing over $k = 1, 2, 3, \dots$ gives a geometric series with first term $\sqrt{3}/12$ and ratio $4/9$. The total area added is
\begin{equation*}
\text{(total area added)} = \frac{\sqrt{3}}{12} \cdot \frac{1}{1 - 4/9} = \frac{\sqrt{3}}{12} \cdot \frac{9}{5} = \frac{3\sqrt{3}}{20}. \tag{15}
\end{equation*}
The total area of the Koch snowflake is therefore
\begin{equation*}
A_\infty = \frac{\sqrt{3}}{4} + \frac{3\sqrt{3}}{20} = \frac{5\sqrt{3}}{20} + \frac{3\sqrt{3}}{20} = \frac{8\sqrt{3}}{20} = \frac{2\sqrt{3}}{5} \approx 0.693. \tag{16}
\end{equation*}
This is the canonical ``infinite perimeter, finite area'' value \cite{Falconer2014,Mandelbrot1982}.

\subsection{Koch-style construction on a square (invented by the author)}\label{sec:koch-square}

The construction presented in this subsection was invented by the author (Pedro Marotta) as part of the original Extended Essay investigation. The iteration rule, the similarity-dimension calculation, the perimeter scaling, and the closed-form area derivation given below are the author's own work.

\begin{figure}[htbp]
  \centering
  \includegraphics[width=0.95\textwidth]{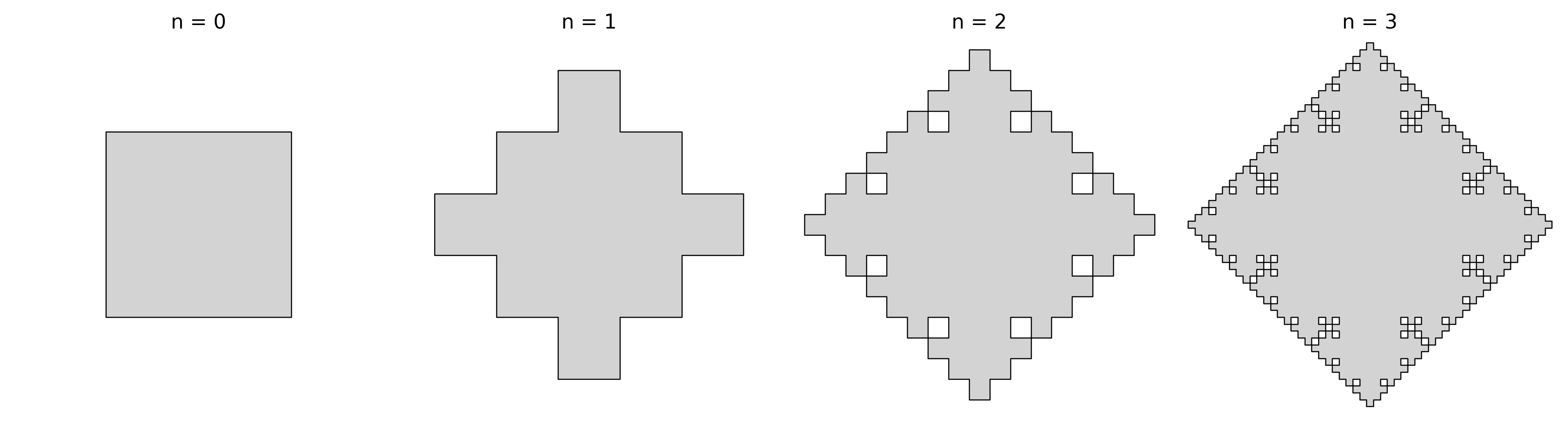}
  \caption{\emph{Koch-style construction on a square, iterations $n=0$ through $n=3$ (author's construction).} The construction starts from a unit square of perimeter~$4$ ($n=0$). At each iteration, every boundary segment of length~$s$ is divided into three equal parts; on the middle part, a square of side $s/3$ is constructed outward; the original middle part is then deleted, so each segment is replaced by five segments of length $s/3$ that form the new outline. The generator parameters are $N=5$, $r=3$, and the similarity dimension of the limit curve is $D = \log 5 / \log 3 \approx 1.465$. The construction is additive: the perimeter diverges as $(5/3)^{n}$, while the area enclosed by the curve converges to $2$ under the explicit non-overlap assumption stated in Section~\ref{sec:koch-square}.}
  \label{fig:koch-square}
\end{figure}

The rule is as follows. The initial figure is a square of side length~$1$, with perimeter~$4$. At each iteration, every segment of the current boundary is divided into three equal parts; on the middle third, a square is constructed outward; the original middle third is then removed (Figure~\ref{fig:koch-square}).

This is an additive construction. Each segment is replaced by five shorter segments (the two flanking thirds, plus three sides of the outward square that constitute the new outline), so $N = 5$; each new segment has one-third the length of its parent, so $r = 3$. The similarity dimension is
\begin{equation*}
D = \frac{\log 5}{\log 3} \approx 1.465. \tag{17}
\end{equation*}
The perimeter at iteration $n$ is
\begin{equation*}
P_n = 4 \cdot \left(\frac{5}{3}\right)^{n}, \tag{18}
\end{equation*}
which diverges as $n \to \infty$.

The area is computed under the assumption that added squares at each iteration intersect neither the existing figure nor each other beyond their bases (the ``non-overlap assumption''; see Section~\ref{sec:nonoverlap} for discussion of this assumption). The initial area is~$1$. At iteration~1, four new squares of side $1/3$ are added (one on each side of the unit square), contributing area $4 \cdot (1/3)^{2} = 4/9$. At iteration~2, every segment of the iteration-1 boundary (twenty segments of length $1/3$) yields one outward square of side $1/9$, so twenty new squares of area $(1/9)^{2} = 1/81$ are added. In general, at iteration $k \ge 1$, $4 \cdot 5^{k-1}$ new squares of side $1/3^{k}$ are added, contributing total area
\begin{equation*}
\text{(area added at iteration } k) = 4 \cdot 5^{k-1} \cdot \left(\frac{1}{3^{k}}\right)^{2} = \frac{4}{9} \cdot \left(\frac{5}{9}\right)^{k-1}. \tag{19}
\end{equation*}
Summing over $k = 1, 2, 3, \dots$ gives a geometric series with first term $4/9$ and ratio $5/9$. The total area added is
\begin{equation*}
\text{(total area added)} = \frac{4}{9} \cdot \frac{1}{1 - 5/9} = \frac{4}{9} \cdot \frac{9}{4} = 1. \tag{20}
\end{equation*}
The total area is therefore
\begin{equation*}
A_\infty = 1 + 1 = 2. \tag{21}
\end{equation*}
The value $A_\infty = 2$ in~(21) is therefore stated as a conditional, not unconditional, result. Under the non-overlap assumption introduced above, the additive computation is exact: the area of the limiting region equals the initial unit-square area $A_0 = 1$ plus the sum of the bump areas given by (19)--(20), which itself sums to~$1$, for a total of~$2$. Without the non-overlap assumption, the same computation provides a finite upper bound rather than the exact value: by countable subadditivity of two-dimensional Lebesgue measure, the area of a union of measurable sets is at most the sum of their individual areas, so the area of the limiting region is at most $A_0$ plus the bump-area series counted with multiplicity, i.e.\ at most $1 + 1 = 2$. The value $A_\infty = 2$ is therefore recorded here as an \emph{unproven yet clearly bounded} claim: unproven, because non-overlap has not been established for all $n$ and so the additive computation cannot be asserted to be exact unconditionally; clearly bounded, both in the structural sense that under non-overlap the closed-form geometric series gives the explicit value~$2$ with no remaining freedom, and in the measure-theoretic sense that~$2$ is in any case a finite upper bound on the limiting area. Section~\ref{sec:nonoverlap} records this contingency and discusses its scope.

\section{The $(N, r)$ framework and regime classification}\label{sec:framework}

\subsection{Scaling laws and derived growth ratios}\label{sec:scaling}

The four examples in Section~\ref{sec:examples} illustrate the general pattern. For a self-similar planar construction with $N$ self-similar pieces per iteration, each scaled by $1/r$, the perimeter at iteration $n$ satisfies
\begin{equation*}
P_n = P_0 \cdot \left(\frac{N}{r}\right)^{n} \tag{22}
\end{equation*}
under the perimeter conventions of \S\ref{sec:perimeter-convention}. The area at iteration $n$ satisfies
\begin{align}
\text{(additive class)} \quad & A_n = A_0 + C \cdot \sum_{k=1}^{n} \left(\frac{N}{r^{2}}\right)^{k-1}, \tag{23} \\
\text{(subtractive class)} \quad & A_n = A_0 \cdot \left(\frac{N}{r^{2}}\right)^{n}, \tag{24}
\end{align}
where $C$ is a positive constant determined by the construction (the contribution of the bumps added at iteration~1 in the additive case).

Two derived growth ratios are introduced as the principal control parameters of the framework:
\begin{align}
\alpha &= \frac{N}{r} \quad \text{(perimeter growth ratio per iteration),} \tag{22a} \\
\beta  &= \frac{N}{r^{2}} \quad \text{(area growth ratio per iteration).} \tag{22b}
\end{align}
In these coordinates, (22) becomes $P_n = P_0 \cdot \alpha^{n}$ and~(24) becomes $A_n = A_0 \cdot \beta^{n}$; the additive area sum~(23) is a geometric series in $\beta$. The asymptotic behavior of the construction is therefore governed jointly by $(\alpha, \beta)$ rather than by similarity dimension alone: $D = \log N / \log r$ coincides with $\log(\alpha r)/\log r = 1 + \log\alpha/\log r$, so $D$ and $\alpha$ are related but distinct, and $\beta$ is not a function of $D$. The framework's organizing observation is that $(\alpha, \beta)$~--- and not $D$ alone~--- are the natural coordinates in which perimeter behavior and area behavior become directly comparable observables.

\subsection{Asymptotic behavior}\label{sec:asymptotic}

The asymptotic behavior of $P_n$ is governed by the perimeter growth ratio $\alpha = N/r$:
\begin{enumerate}
  \item[(i)] $\alpha < 1$ ($N < r$): $P_n$ is bounded (and tends to zero if strict inequality persists at every iteration).
  \item[(ii)] $\alpha = 1$ ($N = r$): $P_n$ is constant equal to $P_0$.
  \item[(iii)] $\alpha > 1$ ($N > r$): $P_n$ diverges to infinity.
\end{enumerate}

The asymptotic behavior of $A_n$ is governed by the area growth ratio $\beta = N/r^{2}$:
\begin{enumerate}
  \item[(i)] $\beta < 1$ ($N < r^{2}$): in the subtractive case $A_n \to 0$; in the additive case the geometric series in~(23) converges and $A_n \to A_0 + C/(1 - \beta)$, a positive finite limit, provided the non-overlap assumption holds.
  \item[(ii)] $\beta = 1$ ($N = r^{2}$): in the subtractive case $A_n$ is constant equal to $A_0$; in the additive case $A_n$ grows linearly in $n$.
  \item[(iii)] $\beta > 1$ ($N > r^{2}$): in the subtractive case $A_n \to \infty$ (which contradicts retention of pieces; this regime is geometrically non-realizable for subtractive constructions of the type considered here); in the additive case $A_n \to \infty$.
\end{enumerate}

The reframing in $(\alpha, \beta)$ coordinates makes the structural point of the framework explicit: a single inequality $\alpha > 1$ governs perimeter divergence, a single inequality $\beta < 1$ governs area boundedness, and the regime classification of \S\ref{sec:regimes} follows from the joint behavior of these two ratios.

\subsection{Three regimes in $(N, r)$-space}\label{sec:regimes}

Combining the perimeter and area asymptotics yields three regimes in the $(N, r)$ parameter space:
\begin{enumerate}
  \item[(a)] \textbf{Subcritical regime ($N \le r$):} perimeter does not diverge. The iterated construction tends to a smooth or weakly irregular limit. Similarity dimension $D \le 1$.
  \item[(b)] \textbf{Intermediate-dimension regime ($r < N < r^{2}$):} perimeter diverges; area behaves differently in the two construction classes, as detailed in~\S\ref{sec:refinement}. Similarity dimension $1 < D < 2$.
  \item[(c)] \textbf{Supercritical regime ($N \ge r^{2}$):} for additive constructions the area grows without bound; for subtractive constructions the construction is not geometrically realizable as a strict subset of the initial figure. Similarity dimension $D \ge 2$.
\end{enumerate}

The condition $r < N < r^{2}$ defining the intermediate-dimension regime is, on taking logarithms, equivalent to $1 < \log N / \log r < 2$, i.e.\ $1 < D < 2$. The intermediate-dimension regime is therefore exactly the set of $(N, r)$ for which the similarity dimension is strictly between~$1$ and~$2$. This condition is well known as the characterization of planar fractal dimension \cite{Falconer2014,Edgar2008} and is presented here as a corollary of the dimension formula~(1) rather than as a new result.

The contribution of the present framework lies not in this corollary but in the change of representation that makes it useful as a joint classification of two observables. Classical analyses of planar self-similar fractals are typically organized around the single observable $D = \log N / \log r$ and the single inequality $1 < D < 2$; the present framework retains those quantities but adds $\alpha = N/r$ and $\beta = N/r^{2}$ as primary coordinates, in which perimeter behavior and area behavior become independently visible. The regime boundaries $N = r$ and $N = r^{2}$~--- equivalently $\alpha = 1$ and $\beta = 1$~--- partition the $(N, r)$ plane into three regions of qualitatively distinct \emph{joint} asymptotic behavior of $P_n$ and $A_n$. The boundaries themselves are not new; their use here as the axes of a joint classification of perimeter and area is the framework's structural point.

\subsection{Construction-class refinement within the intermediate-dimension regime}\label{sec:refinement}

Within the intermediate-dimension regime, the asymptotic two-dimensional Lebesgue measure of the object of interest depends on construction class.

For subtractive constructions, the object of interest is the iterated set itself. By Hutchinson's theorem \cite{Hutchinson1981} on self-similar iterated function systems satisfying the open set condition, any self-similar set with similarity dimension $D < 2$ has two-dimensional Lebesgue measure equal to zero. For both Sierpinski examples, $A_\infty = 0$ follows directly from~(24) since $\beta = N/r^{2} < 1$ in this regime.

For additive constructions, the object of interest is the region bounded by the iterated curve (not the iterated curve itself, which has dimension $D < 2$ and therefore zero two-dimensional Lebesgue measure as a set). Under the non-overlap assumption, the geometric series in~(23) converges, and $A_\infty$ is positive and finite. The Koch snowflake ($A_\infty = 2\sqrt{3}/5$) and the Koch-style construction on the square ($A_\infty = 2$, conditional on the non-overlap assumption recorded in \S\ref{sec:nonoverlap}) illustrate this case.

This is a \emph{structural non-equivalence inside the same dimension class}. Two constructions can have identical $(N, r)$, identical $\alpha$, identical $\beta$, and identical similarity dimension $D$, yet exhibit qualitatively different asymptotic two-dimensional Lebesgue measure depending on which construction class they belong to. Concretely, the Sierpinski triangle ($N = 3$, $r = 2$) and any additive construction with the same parameters would share a similarity dimension of $D = \log 3 / \log 2 \approx 1.585$ yet differ in asymptotic area outcome. The $(N, r, \text{construction class})$ triple therefore carries information that $(N, r)$ alone~--- and equivalently $D$ alone~--- does not. The present framework records this distinction explicitly; it is one of the framework's principal organizing observations.

The non-equivalence is asserted here within the construction classes specified in \S2.2 and is not claimed as a universal property of self-similar planar fractals: constructions that mix additive and subtractive rules, that violate the non-overlap assumption, or that fall outside the deterministic iterated function system class (e.g.\ stochastic constructions; see Section~\ref{sec:discussion}) require separate analysis.

\subsection{Diagnostic use of the framework}\label{sec:diagnostic}

The combined content of \S\S\ref{sec:scaling}--\ref{sec:refinement} may be summarized as a diagnostic procedure. Given a deterministic self-similar planar construction in the class defined in \S2.1, the asymptotic behavior of perimeter and area is determined by the following inputs:
\begin{itemize}
  \item the perimeter growth ratio $\alpha = N/r$, which determines whether $P_n$ is bounded, constant, or divergent (\S\ref{sec:asymptotic});
  \item the area growth ratio $\beta = N/r^{2}$, which determines whether $A_n$ is bounded or divergent in absolute terms (\S\ref{sec:asymptotic});
  \item the construction class (additive or subtractive, \S2.2), which determines, within the intermediate regime, whether the asymptotic two-dimensional Lebesgue measure of the object of interest is zero or positive finite (\S\ref{sec:refinement});
  \item in the additive case, a non-overlap assumption (\S\ref{sec:koch-square}, \S\ref{sec:nonoverlap}) under which the closed-form value of $A_\infty$ is calculable from~(23) as $A_0 + C/(1 - \beta)$.
\end{itemize}

Each of the four worked examples in Section~\ref{sec:examples} is recovered as a special case of this procedure: given $(N, r)$ and construction class, asymptotic perimeter and area follow without re-deriving the geometric series for that example. The framework's contribution is, in this sense, a consolidation of multiple individual derivations into a single diagnostic structure indexed by $(N, r)$ and construction class.

\subsection{Predictive application to new constructions}\label{sec:predictive}

The four examples in Section~\ref{sec:examples} are canonical fractals whose perimeter and area behavior is independently established in the literature, and their treatment in Section~\ref{sec:examples} therefore amounts to a re-derivation of known results in the $(N, r)$ representation rather than to a predictive use of the framework. To demonstrate that the framework also functions predictively~--- i.e.\ that, given the $(N, r, \text{construction class})$ specification of a construction in the class defined in \S2.1--\S2.2, the asymptotic behavior of perimeter and area follows directly from the diagnostic of \S\ref{sec:diagnostic} without performing the full geometric-series derivation case-by-case~--- four further constructions are analyzed below. The four constructions cover the three regimes identified in \S\ref{sec:regimes} (one example per regime) and additionally include a pair with identical $(N, r)$ and identical similarity dimension $D = \log N / \log r$ differing only in construction class, in order to exhibit the structural non-equivalence recorded in \S\ref{sec:refinement}. For each construction the framework's predictions are stated first and a direct calculation is then given to confirm them.

\subsubsection{Subcritical regime: a subtractive $(N, r) = (2, 3)$ construction}

\emph{Construction.} The initial figure is a unit square. At each iteration, every retained sub-square is partitioned into a $3 \times 3$ grid of nine sub-squares of side $1/3$, and only two of the nine are retained~--- for definiteness, the upper-left and lower-right corner sub-squares. The rule is applied recursively to each retained sub-square. The construction is subtractive with $N = 2$ and $r = 3$.

\emph{Diagnostic predictions (from \S\ref{sec:diagnostic}).} The growth ratios are $\alpha = N/r = 2/3 < 1$ and $\beta = N/r^{2} = 2/9 < 1$; the similarity dimension is $D = \log 2 / \log 3 \approx 0.631 < 1$. By \S\ref{sec:asymptotic}, $\alpha < 1$ places the construction in the subcritical regime: the iteration-$n$ perimeter (under the convention of \S\ref{sec:perimeter-convention}) is bounded and tends to zero. By \S\ref{sec:asymptotic}~(i), $\beta < 1$ in the subtractive case yields $A_n \to 0$. The framework therefore predicts that the iteration-$n$ total edge length $P_n$ tends to zero and the iteration-$n$ area $A_n$ tends to zero, consistently with the limit set having two-dimensional Lebesgue measure equal to zero.

\emph{Verification.} Direct substitution into~(22) and~(24) gives $P_n = 4 \cdot (2/3)^{n}$ and $A_n = (2/9)^{n}$, both of which tend to zero as $n \to \infty$. The framework's predictions are confirmed without a separate derivation.

\subsubsection{Intermediate-dimension regime: a new additive $(N, r) = (6, 4)$ construction}

\emph{Construction.} The initial figure is a unit square (perimeter~$4$). At each iteration, every boundary segment of length $\ell$ is divided into four equal parts of length $\ell/4$; on the third part (counted along the segment's direction), an outward square of side $\ell/4$ is constructed; the original third part is then removed. The segment is therefore replaced by six segments of length $\ell/4$ (the first part, the second part, three sides of the outward square, the fourth part). The construction is additive with $N = 6$ and $r = 4$. A non-overlap assumption~--- that the added squares at successive iterations do not intersect either the existing figure or each other beyond their bases~--- is taken to hold, on the same conditional basis recorded in \S\ref{sec:koch-square} and \S\ref{sec:nonoverlap} for the Koch-style construction on a square of \S\ref{sec:koch-square}.

\emph{Diagnostic predictions (from \S\ref{sec:diagnostic}).} The growth ratios are $\alpha = 6/4 = 3/2$ and $\beta = 6/16 = 3/8$; the similarity dimension is $D = \log 6 / \log 4 \approx 1.292$. Since $4 < 6 < 16$, equivalently $r < N < r^{2}$, the construction lies in the intermediate-dimension regime. By \S\ref{sec:asymptotic}, $\alpha > 1$ yields perimeter divergence; $\beta < 1$ in the additive case yields, under the non-overlap assumption, a positive finite asymptotic area $A_\infty = A_0 + C/(1 - \beta)$. The framework therefore predicts simultaneous perimeter divergence and finite positive enclosed area.

\emph{Verification.} From~(22), $P_n = 4 \cdot (3/2)^{n}$, which diverges as $n \to \infty$. For the area, at iteration $k \ge 1$ the construction adds $4 \cdot 6^{k-1}$ outward squares each of area $(1/4^{k})^{2} = 1/16^{k}$, contributing total added area at iteration $k$ equal to $(4/16) \cdot (6/16)^{k-1} = (1/4) \cdot (3/8)^{k-1}$. Summing this geometric series in $\beta = 3/8$ over $k = 1, 2, 3, \dots$ gives total added area $(1/4) \cdot 1/(1 - 3/8) = (1/4) \cdot (8/5) = 2/5$. The asymptotic area is therefore $A_\infty = 1 + 2/5 = 7/5 = 1.4$ under the non-overlap assumption. Both predictions are confirmed.

\subsubsection{Supercritical regime: an additive $(N, r) = (10, 3)$ construction}

\emph{Construction.} The construction is taken to be a formal additive generator with $N = 10$ and $r = 3$: each boundary segment of length $\ell$ is replaced by ten segments of length $\ell/3$, with the bumps assumed to satisfy the same non-overlap convention adopted in \S\ref{sec:koch-square} and \S\ref{sec:nonoverlap} so that the additive area model of~(23) applies. For concreteness one may take the rule in which each segment is divided into three equal parts of length $\ell/3$ and the middle part is replaced by an outward zig-zag of eight segments of length $\ell/3$, so that the segment is replaced by ten segments of length $\ell/3$ in total (the first part, eight zig-zag segments, and the third part); whether such a rule admits a non-self-intersecting planar realization at every iteration is a separate question that the framework, restricted to $(N, r, \text{construction class})$ inputs and the additive area model, does not address.

\emph{Diagnostic predictions (from \S\ref{sec:diagnostic}).} The growth ratios are $\alpha = 10/3 \approx 3.33$ and $\beta = 10/9 \approx 1.11$; the similarity dimension is $D = \log 10 / \log 3 \approx 2.10$. Since $10 > 9 = r^{2}$, the construction lies in the supercritical regime. By \S\ref{sec:asymptotic}, $\alpha > 1$ yields perimeter divergence (under the additive perimeter convention of \S\ref{sec:perimeter-convention}). By \S\ref{sec:asymptotic}~(iii), $\beta > 1$ in the additive case yields divergence of the additive area-counting series~(23), in which bump areas are counted with multiplicity at every iteration. The framework therefore predicts unbounded growth of the iteration-$n$ perimeter $P_n$ and unbounded growth of the additive area-counting series $A_n$.

\emph{Verification.} From~(22), $P_n = 4 \cdot (10/3)^{n}$, which diverges as $n \to \infty$. For the area, the bump-area term at iteration $k$ is proportional to $\beta^{k-1} = (10/9)^{k-1}$, which fails to converge as $n \to \infty$; the additive area sum~(23) therefore diverges. The two predictions made by the framework are confirmed.

The status of this prediction relative to the actual planar Lebesgue measure of the limiting region is the following. The additive series of~(23) counts bump areas with multiplicity and is the natural object on which the framework operates; it is a finite upper bound on the limiting region's area only under the non-overlap assumption (cf.\ \S\ref{sec:koch-square} and \S\ref{sec:nonoverlap}). If non-overlap fails~--- as it generically does in the supercritical regime, where bumps are sufficiently large relative to scale-down to produce overlap or self-intersection~--- then the additive series may diverge while the actual two-dimensional Lebesgue measure of the union of all bump regions remains bounded by countable subadditivity. The framework's prediction is therefore precisely the divergence of the additive area-counting series, not necessarily a statement about the actual Lebesgue measure of any planar realization. The Lebesgue-measure question for supercritical additive constructions falls outside the scope of the $(N, r, \text{construction class})$ framework as stated in \S2.1--\S2.2 and is not addressed here.

\subsubsection{Identical $(N, r)$ and dimension, differing construction class: a $(5, 3)$ pair}

The Koch-style construction on a square presented in \S\ref{sec:koch-square} has parameters $(N, r) = (5, 3)$ and is additive. Consider in addition the following $(5, 3)$ subtractive construction. The initial figure is a unit square. At each iteration, every retained sub-square is partitioned into a $3 \times 3$ grid of nine sub-squares of side $1/3$, and only five of the nine are retained~--- for definiteness, the four corner sub-squares and the central sub-square. The rule is applied recursively. This construction is subtractive with $N = 5$ and $r = 3$.

\emph{Common diagnostic inputs.} Both constructions share $\alpha = 5/3$, $\beta = 5/9$, and $D = \log 5 / \log 3 \approx 1.465$. Both lie in the intermediate-dimension regime $r < N < r^{2}$ of \S\ref{sec:regimes}.

\emph{Diagnostic predictions (from \S\ref{sec:diagnostic}).} The framework predicts identical perimeter behavior for the two constructions: $P_n = P_0 \cdot (5/3)^{n} \to \infty$ in both cases~(22). By the construction-class refinement of \S\ref{sec:refinement}, however, the asymptotic two-dimensional Lebesgue measure differs qualitatively. The additive Koch-style construction on a square gives $A_\infty = 2$ under the non-overlap assumption (\S\ref{sec:koch-square}, \S\ref{sec:nonoverlap}). The subtractive $(5, 3)$ construction satisfies $A_n = (5/9)^{n}$ by~(24), so $A_n \to 0$.

\emph{Verification.} For the additive case, $A_\infty = 2$ has been derived in \S\ref{sec:koch-square} and recorded as a conditional value in \S\ref{sec:nonoverlap}. For the subtractive case, direct calculation: at iteration $n$ there are $5^{n}$ retained sub-squares, each of area $(1/9)^{n}$, so $A_n = (5/9)^{n}$; this tends to zero as $n \to \infty$ since $\beta = 5/9 < 1$. Both predictions are confirmed.

The pair illustrates the structural non-equivalence of \S\ref{sec:refinement} directly. Two constructions with identical $(N, r)$, identical $\alpha$, identical $\beta$, and identical similarity dimension $D$ exhibit qualitatively different asymptotic two-dimensional Lebesgue measures, $A_\infty = 2$ versus $A_\infty = 0$, depending only on construction class. The $(N, r)$ pair, and equivalently $D$ alone, are insufficient to predict the asymptotic area outcome; the $(N, r, \text{construction class})$ triple is required.

\section{Regime figure}\label{sec:regimefigure}

\begin{figure}[htbp]
  \centering
  \includegraphics[width=0.95\textwidth]{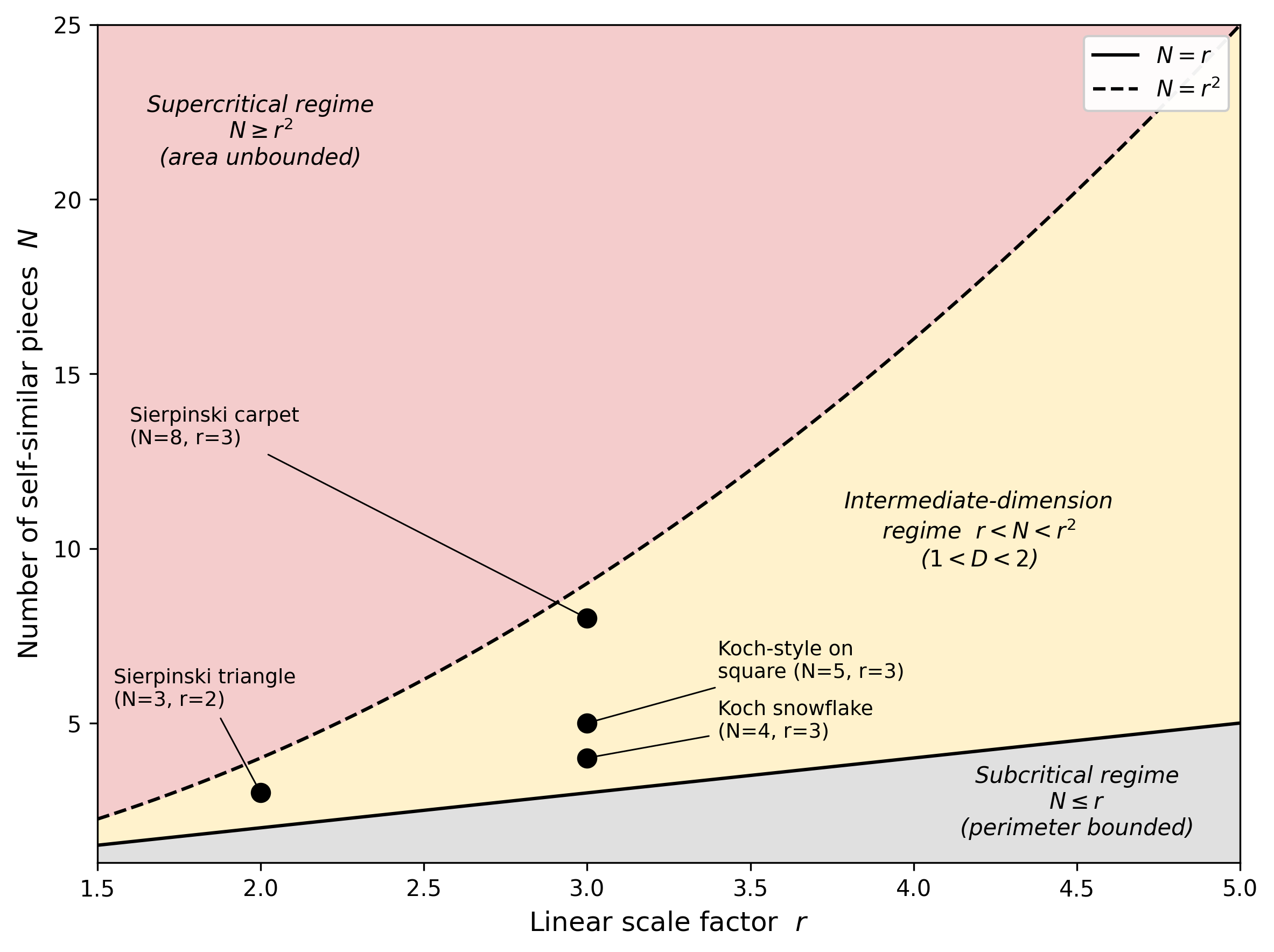}
  \caption{\emph{$(N, r)$ parameter-space regime plot for self-similar planar fractals.} The horizontal axis is the linear scale factor~$r$ and the vertical axis is the number of self-similar pieces~$N$. The two black curves are the boundaries $N = r$ (solid) and $N = r^{2}$ (dashed). The plane is partitioned into three regimes: the subcritical regime $N \le r$ (light grey, similarity dimension $D \le 1$, perimeter bounded under the perimeter convention of Section~\ref{sec:preliminaries}); the intermediate-dimension regime $r < N < r^{2}$ (light yellow, similarity dimension $1 < D < 2$, the regime in which the perimeter convention of this paper diverges and in which the asymptotic two-dimensional Lebesgue measure of the object of interest depends on construction class as discussed in Section~\ref{sec:framework}); and the supercritical regime $N \ge r^{2}$ (light red, similarity dimension $D \ge 2$). The four worked examples are marked: Sierpinski triangle ($N=3$, $r=2$), Sierpinski carpet ($N=8$, $r=3$), Koch snowflake ($N=4$, $r=3$), and the Koch-style construction on a square ($N=5$, $r=3$). All four lie in the intermediate-dimension regime.}
  \label{fig:regime}
\end{figure}

Figure~\ref{fig:regime} plots the $(N, r)$ parameter space with the three regimes shaded and the four worked examples marked: Sierpinski triangle at $(3, 2)$, Sierpinski carpet at $(8, 3)$, Koch snowflake at $(4, 3)$, and the Koch-style construction on the square at $(5, 3)$. The curve $N = r$ marks the boundary between subcritical and intermediate regimes; the curve $N = r^{2}$ marks the boundary between intermediate and supercritical regimes.

The figure is a parameter-space plot, not a phase diagram in the thermodynamic sense, and is presented for visual organization of the regime classification.

\section{Discussion}\label{sec:discussion}

\subsection{Scope of the framework}

The framework presented in this paper applies to the deterministic self-similar iterated function system class with two construction classes (additive and subtractive) as defined in \S2.2. The four worked examples are all special cases of this scheme. The framework is not a universal classification of planar fractals.

\subsection{Stochastic fractal models}

Important fractal phenomena in physics and biology arise from stochastic, not deterministic, processes~--- for example, diffusion-limited aggregation \cite{WittenSander1981} and fractional Brownian motion \cite[Chapter~16]{Falconer2014}. These constructions produce fractal limit sets with non-integer dimensions but are governed by random rather than deterministic iteration rules. The $(N, r)$ parameterization adopted here does not apply to these models, and any apparent visual resemblance between deterministic IFS fractals and stochastic fractals does not imply a shared mathematical framework. This paper makes no claims about stochastic fractal models.

\subsection{Non-overlap assumption}\label{sec:nonoverlap}

The asymptotic area calculation for additive constructions (e.g.\ the Koch-style construction on the square in \S\ref{sec:koch-square}) relies on a non-overlap assumption: that bumps added at each iteration intersect neither the existing iterated region nor each other beyond their bases. For the Koch snowflake, non-overlap follows from the geometry of the equilateral triangle and the inward-pointing angles of the construction. For the Koch-style construction on the square, non-overlap holds for the first three iterations as verified visually in Figure~\ref{fig:koch-square}; a general proof for arbitrary iteration depth requires checking that bumps added on perpendicular boundary segments at successive iterations do not collide.

The status of the value $A_\infty = 2$ in~(21) is therefore explicitly the following. It is an \emph{unproven yet clearly bounded} claim. \emph{Unproven} in the sense that the non-overlap assumption on which the geometric-series derivation depends has been verified visually for iterations $n = 0, 1, 2, 3$ only, and a general argument for arbitrary $n$ is not provided in this paper; visual verification at finite $n$ is consistent with non-overlap holding for all $n$ but does not establish it. \emph{Clearly bounded} in two senses. First, structurally: conditional on non-overlap, the additive computation $A_0 + (\text{sum of bump areas})$ yields the explicit value $1 + 1 = 2$ with no remaining freedom~--- there is no unspecified constant, no asymptotic equivalence, and no parameter that can shift the value. Second, measure-theoretically: even if non-overlap fails, two-dimensional Lebesgue measure is countably subadditive, so the area of the limiting region is at most the initial-square area plus the sum of the bump areas counted with multiplicity. The bump-area series in~(19)--(20) sums to~$1$, and adding the initial area $A_0 = 1$ gives a total of~$2$; the value~$2$ is therefore a finite upper bound on the actual limiting area regardless of whether the bumps overlap. The two readings agree under non-overlap and bracket the value otherwise. The value $A_\infty = 2$ is therefore stated throughout this paper as a closed-form conditional result, exact under the stated non-overlap assumption and a finite upper bound on the limiting area without it, rather than as an unconditional theorem. Establishing or refuting non-overlap for the construction of \S\ref{sec:koch-square} at arbitrary iteration depth is left as future work.

\subsection{Perimeter convention for subtractive constructions}

The convention adopted in \S\ref{sec:perimeter-convention} for the perimeter of subtractive constructions~--- total edge length of all sub-pieces present at iteration $n$~--- is non-standard. Alternative conventions yield different values: the boundary of the convex hull of the iterated set is constant; the Hausdorff one-measure of the limit set is infinite for sets of dimension greater than one; the topological boundary of the limit set differs from both. The ``diverging perimeter'' claimed for the Sierpinski examples in \S\ref{sec:sierpinski-triangle} and \S\ref{sec:sierpinski-carpet} is specifically with respect to the total-edge-length convention, and a different convention would yield a different scaling.

\subsection{Hausdorff vs similarity dimension}

The similarity dimension $D = \log N / \log r$ used throughout this paper coincides with the Hausdorff dimension of the limit set for self-similar iterated function systems satisfying the open set condition \cite{Hutchinson1981}. The general definition of Hausdorff dimension is measure-theoretic \cite[Chapter~2]{Falconer2014} and is in general distinct from the similarity dimension. For all four worked examples, the open set condition is satisfied and the two dimensions coincide; the calculations reported here are similarity-dimension calculations.

\section{Conclusion}\label{sec:conclusion}

This paper has presented a unified parameter-space representation for a class of deterministic self-similar planar constructions, organized by the pair $(N, r)$ and the derived growth ratios $\alpha = N/r$ and $\beta = N/r^{2}$. Within this representation, perimeter behavior and area behavior become independently visible as observables governed by $\alpha$ and $\beta$ respectively, and the asymptotic behavior of any construction in the class can be determined directly from its parameters and construction class without re-deriving the geometric series for that example.

The condition $r < N < r^{2}$~--- equivalent to similarity dimension $1 < D < 2$~--- defines the intermediate-dimension regime in which the perimeter convention adopted here diverges. Within that regime, the framework adds a construction-class refinement: under the non-overlap assumption (\S\ref{sec:koch-square}, \S\ref{sec:nonoverlap}), additive constructions yield positive finite asymptotic two-dimensional Lebesgue measure, while subtractive constructions yield zero. This refinement records a structural non-equivalence inside the same dimension class~--- two constructions with identical $(N, r)$, and therefore identical $D$, can exhibit qualitatively different asymptotic area outcomes depending on construction class~--- that is not visible from $D$ alone (\S\ref{sec:refinement}).

The contribution of this paper is one of formulation and synthesis rather than of new mathematics. The dimension formula $D = \log N / \log r$ and the condition $1 < D < 2$ are classical \cite{Falconer2014,Edgar2008}; Hutchinson's theorem on the two-dimensional Lebesgue measure of self-similar sets with $D < 2$ is standard \cite{Hutchinson1981}. What is offered here is (i)~a \emph{change of representation} in which $(\alpha, \beta)$ are the primary coordinates and the perimeter and area observables become directly comparable, (ii)~a \emph{synthesis} in which four canonical examples~--- the Sierpinski triangle, the Sierpinski carpet, the Koch snowflake, and a Koch-style construction on a square invented by the author~--- are located within a single coordinate system rather than analyzed in isolation, and (iii)~a \emph{construction-class refinement} that records the structural non-equivalence described in \S\ref{sec:refinement}. The framework is intended to be used as a diagnostic (\S\ref{sec:diagnostic}): given $(N, r)$ and construction class, perimeter divergence, area boundedness, and the zero-versus-positive-finite character of the asymptotic Lebesgue measure can be inferred directly, subject to the stated assumptions.

The diagnostic use of the framework is exhibited in \S\ref{sec:predictive} on four further constructions outside the canonical four of Section~\ref{sec:examples}: a subtractive $(2, 3)$ construction in the subcritical regime, a new additive $(6, 4)$ Koch-style construction in the intermediate regime (giving $A_\infty = 7/5$ under non-overlap), an additive $(10, 3)$ construction in the supercritical regime (perimeter and the additive area-counting series both unbounded), and a $(5, 3)$ pair~--- the additive Koch-style construction on a square of \S\ref{sec:koch-square} ($A_\infty = 2$ under non-overlap) and a subtractive $(5, 3)$ construction ($A_\infty = 0$)~--- that share $\alpha$, $\beta$, and similarity dimension $D$ yet differ qualitatively in asymptotic Lebesgue measure. In each of these four cases the framework's predictions follow from $(N, r, \text{construction class})$ and are subsequently confirmed by direct calculation, so that the framework is exercised predictively and not only post hoc.

The Koch-style construction on a square presented in \S\ref{sec:koch-square} illustrates the additive-class case at $(N, r) = (5, 3)$, giving the closed-form value $A_\infty = 2$ under the non-overlap assumption. The status of this value is recorded explicitly in \S\ref{sec:koch-square} and \S\ref{sec:nonoverlap} as an unproven yet clearly bounded conditional claim: closed-form and numerically explicit, exact under the stated non-overlap assumption (verified visually for low iterations but not proved for arbitrary iteration depth) and, by countable subadditivity of two-dimensional Lebesgue measure, a finite upper bound on the limiting area without that assumption.

Returning to the question that motivated the original investigation~--- when does a self-similar planar construction yield an iterated curve of divergent perimeter bounding a region of positive finite area~--- the framework provides an explicit diagnostic answer: an additive construction in the intermediate-dimension regime $r < N < r^{2}$ with the non-overlap assumption holding. The Koch snowflake and the Koch-style construction on the square are two such examples; the Sierpinski triangle and Sierpinski carpet, although in the same regime, fall in the subtractive class and have asymptotic area zero. The scope of the framework is restricted to deterministic self-similar constructions in the two classes defined in \S2.2; mixed constructions, overlapping geometries, and stochastic fractal models lie outside its scope (\S\S\ref{sec:discussion}.1--\ref{sec:discussion}.2) and are not addressed here.

\section*{Acknowledgments}

Reviewer~1 of the original submission of this paper to the \emph{Journal of High School Science} suggested the $(N, r)$ parameterization that is adopted in this revision as the paper's organizing framework. The same reviewer, in a second round of comments, suggested reframing the contribution as a representation-and-synthesis result and presenting the framework as a diagnostic tool; that reframing is adopted in the round-2 revision and is reflected in the Abstract, \S1, \S\ref{sec:scaling}, \S\ref{sec:diagnostic}, and \S\ref{sec:conclusion}. AI tools (ChatGPT and Claude) were used during preparation of these revisions for editorial assistance, structural suggestions, and checking of algebraic exposition. All mathematical derivations and the Koch-style construction on a square presented in Section~\ref{sec:koch-square} are the author's own work. The author verified all mathematical content and accepts full responsibility for the manuscript.

\section*{Code availability}

Python code for generating the iterations of the Sierpinski triangle, Sierpinski carpet, and Koch snowflake (originally provided as Appendices A--C in the Journal of High School Science version of this paper), together with the matplotlib scripts that generated Figures~1--5, is openly available at \url{https://github.com/pedromarotta/fractal-framework-demo}.

\end{document}